\documentclass[onecolumn]{article}
\usepackage{nolta2014}
\usepackage{txfonts}

\usepackage{graphicx}

\def\rz{\ifmmode {I\hskip -2.1pt R} \else {\hbox {$I\hskip -2pt R$}}\fi}

\def\cz{\ifmmode {C\hskip -4.5pt\vrule height5.0pt\hskip 4.3pt} \else
       {\hbox {$C\hskip -4.5pt\vrule height5.0pt\hskip 4.3pt$}}\fi}

\usepackage[utf8]{inputenc}         
\usepackage{pict2e}               

\def\ODM{ 
\mbox{
\setlength{\unitlength}{7pt}
\begin{picture}(1,1)
\put(1,1){\line(-1,0){1}}
\put(1,1){\line(0,-1){1}}
\put(0,1){\line(1,-1){1}}
\end{picture}
}}

\usepackage{array}
\usepackage{booktabs}
\usepackage{multirow}
\usepackage{bigdelim}

\newcolumntype{_}{>{\global\let\currentrowstyle\relax}}
\newcolumntype{^}{>{\currentrowstyle}}

\begin{document}

\title{A Contribution to the Numerics of Polynomials and Matrix Polynomials}

\author{Sigurd Falk ${}^\dag$}

\address{
\dag  Leibniz Universit\"at Hannover, Germany \\ 
Institute of Theoretical Electrical Engineering \\
Email: falk@tet.uni-hannover.de
}

\maketitle

\begin{abstract}
In this paper some algorithms will be presented which can be used for the calculation of zeros of polynomials and eigenvalues of polynomial matrices with a multiplicity larger than one. The numerical values calculated with MATLAB are used as starting values. The reliability of the algorithms is demonstrated by means of 8 examples. 
\end{abstract}

\section{Formulation of the Problem}

\noindent Let be a matrix eigenvalue equation 
\begin{equation}\label{MEE}
\mathbf{y}^T \mathbf{F}(\lambda)=\mathbf{0}^T,\quad \mathbf{F}(\lambda) \mathbf{x} =\mathbf{0}
\end{equation}
with a polynomial matrix of the order $n$ and the degree $\rho$
\begin{equation}\label{PM}
\mathbf{F}(\lambda) = \mathbf{A}_0 + \mathbf{A}_1 \lambda + \mathbf{A}_2 \lambda^2 + \cdots + \mathbf{A}_\rho \lambda^\rho;\quad \det \mathbf{A}_\rho \not= 0
\end{equation}
and complex-valued coefficient matrices $\mathbf{A}_0,\ldots , \mathbf{A}_\rho$. 
\vspace{0.5cm}

\noindent In the following the eigenvalues of $\mathbf{F}(\lambda)$
\begin{equation}
\lambda_1,\lambda_2,\ldots , \lambda_m; \quad m=\rho\cdot n
\end{equation}
defined as zeros of the characteristic polynomial
\begin{equation}\label{determinant}
\det \mathbf{F}(\lambda) = f(\lambda)=a_0 + a_1 \lambda + a_2 \lambda^2 + \cdots + a_m \lambda^m
\end{equation}
will be calculated by means of a Pad\'e function\footnote{Henri Eugene Pad\'e, French mathematician, 1863-1953}
\begin{equation}
p(\lambda)=\frac{f(\lambda)}{z(\lambda)},
\end{equation}
where $z(\lambda)$ is a polynomial of degree $\leq m$. 
\begin{figure}[h]
\centering
\includegraphics[width=0.7\linewidth]{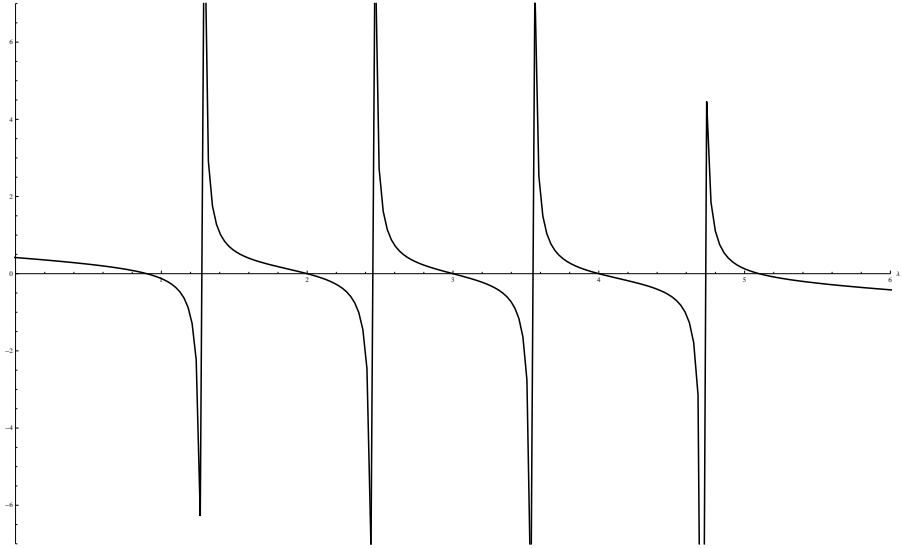}
\caption[]{Typical graph of a Pad\'e function}
\label{fig:Falk}
\end{figure}

\noindent Choosing $z(\lambda)=-f'(\lambda)$ the Pad\'e function
\begin{equation}\label{Pade}
p(\lambda)=\frac{f(\lambda)}{-f'(\lambda)}
\end{equation}
with an interesting property is obtained.

\noindent If the polynomial $f(\lambda)$ possesses a zero $a$ with the multiplicity $\nu$ then it can be represented by 
\begin{equation}
f(\lambda)= (\lambda -a)^\nu\cdot z(\lambda),\quad z(a)\not= 0.
\end{equation}
It follows 
\begin{equation}
f'(\lambda)= \nu (\lambda -a)^{\nu -1}\cdot z(\lambda) +(\lambda -a)^\nu\cdot z'(\lambda),
\end{equation}
or 
\begin{equation}
f'(\lambda)= (\lambda -a)^{\nu -1} \left[ \underbrace{\nu\cdot z(\lambda)+ (\lambda -a)\cdot z'(\lambda) }_{\varphi(\lambda )}      \right] ,
\end{equation}
where 
\begin{equation}
\varphi (a)= \nu\cdot z(a)+0\not= 0 .
\end{equation}

\noindent Therefore, the following theorem can be formulated:
\vspace{0.3cm}

\noindent The Pad\'e function (\ref{Pade}) possesses only zeros with the multiplicity $\nu = 1$. 

\section{Algorithms}\label{algorithms}

\noindent For users only such algorithms are of interest which calculate zeros also if their multiplicity is larger than one. Therefore, Newton's tangent method has to be excluded. In numerical applications there are three methods that can be used in a successful manner. \vspace{0.5cm}

\noindent 2a) A method that is based on the Pad\'e function (\ref{Pade}) leads to
\begin{equation}\label{algor21}
\Lambda_{j+1}=\Lambda_{j}+p(\Lambda_{j}); \quad j=1,2,\ldots
\end{equation}
\noindent 2b) A method that is founded on Halley's function
\begin{equation}
h(\lambda)=\frac{p(\lambda)}{1+p(\lambda)\cdot q(\lambda)}
\end{equation} 
where
\begin{equation}
q(\lambda)=\frac{f''(\lambda)}{f'(\lambda)}.
\end{equation}
leads to
\begin{equation}\label{Halley}
\Lambda_{j+1}=\Lambda_{j}+h(\Lambda_{j}); \quad j=1,2,\ldots
\end{equation}
\noindent 2c) A method that is based on the test polynomials 
\begin{equation}
f_k(\lambda)= a_0 + (-1)^k\cdot a_2\,\lambda^2 + (-2)^k\cdot a_3\,\lambda^3 +\cdots + 
(1-m)^k\cdot a_m\,\lambda^m;
\end{equation}
$k=1,2,\ldots , \nu $, where $\nu$ is the multiplicity of the zero under consideration. \vspace{0.5cm}

\noindent If the prescribed polynomial $f_0(\lambda)$ possesses a zero $\tilde \lambda$ with the multiplicity $\nu$ then each of the Pad\'e functions
\begin{equation}
P_1(\lambda)=\frac{f_0(\lambda)}{f_1(\lambda)},\ P_2(\lambda)=\frac{f_1(\lambda)}{f_2(\lambda)},\ldots , P_\nu(\lambda)=\frac{f_{\nu - 1}(\lambda)}{f_\nu(\lambda)}
\end{equation} 
possesses this zero with the multiplicity one. With
\begin{equation}
p_\nu(\lambda)= P_\nu(\lambda)\cdot \lambda
\end{equation}
the algorithm can be formulated by
\begin{equation}\label{Lambdapv}
\Lambda_{j+1}=\Lambda_{j}+p_\nu(\Lambda_{j}); \quad j=1,2,\ldots
\end{equation}
or 
\begin{equation}\label{Lambda_1+pv}
\Lambda_{j+1}=\left[ 1 +P_\nu(\Lambda_{j})\right]\cdot \Lambda_{j} ; \quad j=1,2,\ldots
\end{equation}
All three algorithms convergence quadratic if the starting value is chosen in a suitable interval that includes the desired zero. This condition is fulfilled if MATLAB results are used as starting values. 

\section{Exploration}\label{exploration}

\noindent In order to obtain suitable approximated values for the start of the algorithms an exploration is needed where three cases 3a1), 3a2) and 3b) have to be distinguish. \vspace{0.5cm}

\noindent 3a) The coefficients of the prescribed polynomial (\ref{determinant}) as well as the zeros are real. Then, the Pad\'e function (\ref{Pade}) is used. 

\noindent 3 a1) The usual regula falsi method.

\noindent If for two arbitrary test points $\lambda_1$ and $\lambda_2$ a change of sign occur with
\begin{equation}\label{sign}
\lambda_1<\lambda_2; \quad p(\lambda_1)>0, p(\lambda_2)<0 ,
\end{equation}
then a zero of the function $p(\lambda)$ is placed between $\lambda_1$ and $\lambda_2$ and therefore also a zero of $f(\lambda)$ exists possibly with a multiplicity larger than one. 

\noindent Using the regula falsi method 
\begin{equation}\label{delta2}
\lambda_3=\lambda_1-\frac{p(\lambda_1)}{\Delta_2}
\end{equation}
with the difference quotient 
\begin{equation}
\Delta_2 = \frac{p(\lambda_2) - p(\lambda_1)}{\lambda_2 - \lambda_1}
\end{equation}
a first approximated  value (in general crude) for a zero is received.
\vspace*{0.5cm}

\noindent 3 a2) Regula falsi method with acceleration

\noindent The approximated value $\lambda_3$ can be improved in the following manner.  If a further difference quotient 
\begin{equation}\label{delta3}
\Delta_3 = \frac{p(\lambda_3) - p(\lambda_1)}{\lambda_3 - \lambda_1}
\end{equation}
and the terms 
\begin{equation}\label{Q2}
Q_2=\frac{p(\lambda_2)}{p(\lambda_1)}, \ Q_3=\frac{p(\lambda_3)}{p(\lambda_1)}
\end{equation}
are defined an improved approximated value $\lambda_4$ is calculated by 
\begin{equation}\label{lambda4}
\lambda_4 = \lambda_1 -\frac{p(\lambda_2) - p(\lambda_3)}{Q_2 \Delta_3 - Q_3 \Delta_2}
\end{equation}
and this scheme can be continued as follows
\vspace{0.5cm}

\noindent  First Step. 

Replace in (\ref{delta2}) to (\ref{lambda4}) the indexes 1,2 and 3 through 2,3 and 4 and calculate $\lambda_5;\, p(\lambda_5)$.
\vspace{0.2cm}

\noindent  Second Step.

Replace in (\ref{delta2}) to (\ref{lambda4}) the indexes 1,2 and 3 through 3,4 and 5 and calculate $\lambda_6;\, p(\lambda_6)$.
\vspace{0.5cm}

\noindent This iteration process has to be broken if 

\noindent a) the condition
\begin{equation}\label{psigma}
\vert p(\lambda_\mu)\vert \leq 10^{-\sigma}
\end{equation}

is fulfilled or 

\noindent b) stop the iteration at a certain iteration step $\mu$ without considering a stopping criteria.
\vspace{0.5cm}

\noindent Now, we discuss the exploration process. After choosing a step-size $\delta$ we calculate on the $\lambda$-axis pairs of values 
\begin{equation}
\lambda_j;\, p(\lambda_j);\quad j=1,2,\ldots
\end{equation}
beginning from zero until a first, second, third, etc. change of sign is found. 
\vspace{0.5cm}

\noindent 3 b) In order to calculate also the negative zeros the co-function 
\begin{equation}\label{PadeCo}
\hat p(\lambda)=-\frac{f(-\lambda)}{-f'(-\lambda)}
\end{equation}
the sequence of steps (\ref{delta3}). (\ref{Q2}) and (\ref{lambda4}) have to be carried out until all $m$ zeros are calculated. 
\vspace{0.5cm}

\noindent 3c) Much more tedious is the exploration within the complex plane since no change of sign in the sense of (\ref{sign}) is available. 
\vspace{0.5cm}

\noindent 3d) Diagonal dominant polynomial matrices

\noindent In the case of distinct diagonal dominance of a matrix the $m=\rho\cdot n$ zeros of the equations
\begin{equation}\label{quadraticEq}
f_{jj}(\lambda)=0; \quad j=1,2,\ldots ,n
\end{equation}
are suitable starting points for the in section 2 presented algorithms. In the case of $\rho =2$ we have to solve $n$ quadratic equations; see also the second example. 

\section{Eigenvalues of a Polynomial Matrix}

\noindent In the following we consider equation (\ref{MEE})
\begin{equation}
\mathbf{F}(\lambda) \mathbf{x} =\mathbf{0}
\end{equation}
with the polynomial matrix (\ref{PM}). 

\noindent Let be $\lambda_k$ an eigenvalue with the multiplicity one, then the matrix 
\begin{equation}\label{Flambda_k}
\mathbf{F}(\lambda_k)
\end{equation}
has the rank $n-1$. 
\vspace{0.5cm}

\noindent a) Transformation of Gau\ss

\noindent If necessary a column pivot search as well as the changing of two rows will be arranged such that the matrix (\ref{Flambda_k}) has the form 
\begin{equation}\label{Ftriangle}
\tilde\mathbf{F}(\lambda_k)= 
 \pmatrix{\tilde{\ODM}_k & \mathbf{w}_k\cr
 \mathbf{0}^T & 0 }
\end{equation}
where  $\tilde{\ODM}_k$ is  regular upper rectangular matrix of the order $n-1$ and the column $\mathbf{w}_k$ has the length $n-1$.
\vspace*{0.5cm}

\noindent b) Transformation due to Jordan in the order\footnote{Wilhelm Jordan, Geometer, 1842-1899}
\begin{equation}
n-1,n-2,\ldots, 2.
\end{equation}
Therefore, 
\begin{equation}
\hat{\mathbf{F}}(\lambda_k)= \pmatrix{\mathbf{D}_k & \mathbf{z}_k \cr \mathbf{0}^T& 0}
\end{equation}
with a regular diagonal matrix $\mathbf{D}_k$ of the order $n-1$.
\vspace{0.2cm}

\noindent It is easy to see that the desired eigenvector is 
\begin{equation}\label{EV_Ftilde}
\mathbf{x}_k=\pmatrix{\mathbf{D}_k^{-1} \mathbf{z}_k\cr -1}.
\end{equation}
Moreover, since the system of equations (\ref{MEE}) is homogeneous, 
\begin{equation}\label{EVnorm1}
\hat{\mathbf{x}}_k=\alpha_k\cdot \mathbf{x}_k,\quad \alpha_k\not= 0
\end{equation}
is also an eigenvector. The factor $\alpha_k$ can be determined such that $\hat{\mathbf{x}}_k$ is orthonormal 
\begin{equation}\label{EVnorm2}
\hat{\mathbf{x}}_k^*\hat{\mathbf{x}}_k =1,
\end{equation}
but we choose 
\begin{equation}\label{EVnorm3}
\alpha_k=1.
\end{equation}
\vspace{0.5cm}

\noindent Multiple Eigenvalues

\noindent Let be $\lambda_k$ an eigenvalue with the multiplicity $\nu_k$ and $r_k$ the rank deficiency of the matrix $\mathbf{F}(\lambda_k)$, where we have
\begin{equation}
r_k\leq \nu_k .
\end{equation}
Now, the matrix (\ref{Ftriangle}) has the form
\begin{equation}
\tilde\mathbf{F}(\lambda_k)=\pmatrix{\tilde{\ODM} & \mathbf{0}\cr
 \mathbf{0}^T & \mathbf{0} }
\end{equation}
where the zero matrix in the right lower corner has the order $r_k$. 

\noindent Corresponding to (\ref{Ftriangle}) to (\ref{EV_Ftilde}) we have 
\begin{equation}
\tilde\mathbf{F}(\lambda_k)=\pmatrix{\mathbf{D}_k& \tilde{\mathbf{Z}}_k\cr
\mathbf{0} & \mathbf{0}}
\end{equation}
and therefore
\begin{equation}\label{EVrk}
\mathbf{X}_k =\pmatrix{\mathbf{D}_k^{-1} \tilde{\mathbf{Z}}_k\cr -\mathbf{I}_{r_k}}
=\left( \mathbf{x}_1\ \mathbf{x}_2\ \cdots \mathbf{x}_{r_k} \right) .
\end{equation}
The vectors $\mathbf{x}_i\ (i=1,\ldots ,r_k)$ are the $r_k$ linear independent eigenvectors of $\lambda_k$ that can be normed with respect to (\ref{EVnorm1}) - (\ref{EVnorm3}). 
\vspace{0.2cm}

\noindent If we have $r_k< \nu_k$ then the $r_k$ eigenvectors (\ref{EVrk}) can be complemented by generalized eigenvectors; cf. \cite{ZurmuehlFalk2}.
\vspace{0.5cm}

\noindent Now, we consider the left eigenvectors
\begin{equation}
\mathbf{y}^T\, \mathbf{F}(\lambda)=\mathbf{0}
\end{equation}
where after a transposition of this equation it follows
\begin{equation}
\left[\mathbf{y}^T\, \mathbf{F}(\lambda)\right]^T =\left[\mathbf{0}^T\right]^T \ \Rightarrow\ \mathbf{F}(\lambda)^T \mathbf{y}=\mathbf{0}
\end{equation}
If $\mathbf{F}$ is replaced by $\mathbf{F}^T$ the concepts of this section can be used. 

\section{The EPC-Transformation}

\noindent The algorithm described in \cite{ZurmuehlFalk2} based on the allocation of $m$ pairwise different interpolation values
\begin{equation}\label{interpvalues}
\sigma_1,\sigma_2,\ldots , \sigma_m,
\end{equation}
which have to be chosen in suitable manner. With these values the following interpolation polynomials are defined
\begin{equation}
g_k(\lambda)=\prod_{j=1 \atop j\not= k }^m (\sigma_j - \lambda )  ;\quad k=1,2,\ldots , m
\end{equation}
and therefore the Pad\'e functions
\begin{equation}
P_k(\lambda)=\frac{f(\lambda)}{g_k(\lambda)}\cdot \frac{1}{a_m};\quad k=1,2,\ldots , m .
\end{equation}
For $\lambda =\sigma_k$ we obtain the defects (as denoted in \cite{CarstensenStein1987} and \cite{CarstensenStein1989})
\begin{equation}\label{defects}
d_k=\frac{f(\sigma_k)}{g_k(\sigma_k)}\cdot \frac{1}{a_m};\quad k=1,2,\ldots , m .
\end{equation}
and the corresponding so-called main values
\begin{equation}
H_k=\sigma_k-d_k;\quad k=1,2,\ldots , m .
\end{equation}
These values will be collected in the following list
\begin{equation}\label{list}
L_m=\pmatrix{\hbox{Interpolation Values}&\hbox{Defects}&\hbox{Main Values}\cr 
\midrule\cr
\sigma_1&d_1& H_1\cr
\sigma_1&d_1& H_1\cr
\vdots & \vdots & \vdots \cr
\sigma_m&d_m& H_m},
\end{equation}
that includes the entire information of the polynomial matrix (\ref{PM}). The order of
the rows is arbitrary. The control equation
\begin{equation}\label{SumHj}
\sum_{j=1}^m H_j = -\frac{a_{m-1}}{a_m}
\end{equation}
error-free calculation of the defects (\ref{defects}) from the interpolation values (\ref{interpvalues}). 

\section{The ECP-Rayleigh Quotient}

\noindent Let be the eigenvalue equation (\cite{Falk2004}, p. 421)
\begin{equation}\label{EVE}
\det \mathbf{F}(\lambda) = \det (\mathbf{E} - \lambda \mathbf{I}_m)=0
\end{equation}
with the accompanying ECP matrix
\begin{equation}
\mathbf{E} = \hbox{Diag}< \sigma_j> - \pmatrix{1\cr 1\cr \vdots\cr 1} (d_1\ d_2\ \cdots \ d_m).
\end{equation}
The Rayleigh quotient 
\begin{equation}
R(\lambda)=\frac{\mathbf{y}^T(\lambda) \mathbf{E}\, \mathbf{x}(\lambda)}{\mathbf{y}^T(\lambda) \mathbf{I}_m\, \mathbf{x}(\lambda)},
\end{equation}
where
\begin{eqnarray}
\mathbf{y}^T(\lambda)&=&\pmatrix{ \frac{1}{\sigma_1-\lambda}& \frac{1}{\sigma_2-\lambda}& \cdots & \frac{1}{\sigma_m-\lambda}},\\
\\
\mathbf{x}(\lambda)&=& 
\pmatrix{
\frac{d_1}{\sigma_1 -\lambda }\cr 
\frac{d_2}{\sigma_2 -\lambda }\cr 
\vdots \cr 
\frac{d_m}{\sigma_m-\lambda }
}
\end{eqnarray}

\noindent can be reformulated by using the terms 
\begin{eqnarray}
S_1(\lambda ) &=&  \sum_{j=1}^m \frac{d_j}{\sigma_j-\lambda} , \\
S_2(\lambda ) &=&  \sum_{j=1}^m \frac{d_j}{(\sigma_j-\lambda )^2} , \\
S_\sigma(\lambda ) &=&  \sum_{j=1}^m \frac{d_j \sigma_j}{(\sigma_j-\lambda )^2}   
\end{eqnarray}
in the form
\begin{equation}
R(\lambda)= \frac{S_\sigma(\lambda) - S_1^2(\lambda)}{S_2(\lambda)} .
\end{equation}
Therefore, the following algorithm is defined
\begin{equation}
\Lambda_{j=1}=\Lambda_j+R(\Lambda_j); \quad j=1,2,\ldots 
\end{equation}
which can be started by a main value $H_k$. 

\section{The Reduced Eigenvalue Equation}

\noindent Among the eigenvalue equation (\ref{EVE}) the reduced eigenvalue equation exists according to (\cite{Falk2004}, p. 346)
\begin{equation}\label{detF}
\tilde f (\lambda ) = S_1 (\lambda ) -1 =0.
\end{equation}
With the derivative 
\begin{equation}
\tilde f' (\lambda ) = S'_1 (\lambda ) -0 =S_2 (\lambda )
\end{equation}
the Pad\'e function
\begin{equation}
p_E(\lambda)=\frac{\tilde f (\lambda )}{-\tilde f' (\lambda )} =\frac{ S_1 (\lambda ) -1}{-S_2 (\lambda )} 
\end{equation}
is obtained and therefore the algorithm 
\begin{equation}\label{redEV}
\Lambda_{j=1}=\Lambda_j+p_E(\Lambda_j); \quad j=1,2,\ldots 
\end{equation}
It can be started by a main value $H_k$. 

\section{The Evolution}\label{evolution}

\noindent  8a) An additional algorithm is introduced in (\cite{Falk2004}, p.44) which can be described as follows: replace the interpolation values in list (\ref{list}) by the main values and prepare a new list; repeat this procedure as long as some or all defects go below a prescribed threshold. The main values of the final list can be used as start values of the algorithm in section \ref{algorithms}. 

\section{Numerical Feasibility and Additional Aspects}

\noindent 9a) Evaluation of the multiplicity for the algorithm (\ref{Lambdapv}). 
\vspace{0.5cm}

\noindent Execute the algorithm for $\nu =1$, $\nu =2$, and so on, simultaneously. The convergence will be taken place exactly once. Therefore, the zeros and their multiplicity is determined. 

\noindent The Taylor test with the characteristic polynomial (\ref{determinant})
\begin{eqnarray}
\matrix{\nu =1\colon & f(a) & =& 0\cr
&f'(a)&\not= & 0\cr
\nu =2\colon & f(a) & =& 0\cr
&f'(a)&= & 0\cr
&f''(a)&\not= & 0\cr}
\end{eqnarray} 
and the same manner for $\nu>2$ can be used as control.
\vspace{0.5cm}

\noindent 9b) The matrix (\ref{EVE}) can be reformulated as 
\begin{equation}
\mathbf{E} = 
\pmatrix{
H_1& -d_2& -d_3 & \cdots& -d_{m-1} & -d_m \cr
-d_1& H_2& -d_3 & \cdots& -d_{m-1} & -d_m \cr
-d_1 & -d_2 & H_3& \cdots& -d_{m-1} & -d_m \cr
\vdots & \vdots & \vdots & \ddots & \vdots &\vdots \cr
-d_1 & -d_2 & -d_3& \cdots& H_{m-1} & -d_m \cr
-d_1 & -d_2 & - d_3& \cdots& -d_{m-1} & H_m \cr
}
\end{equation}
and therefore 
\begin{equation}
\mathrm{Tr}\hspace{2pt} \mathbf{E} = \sum_{j=1}^m H_j=  \sum_{j=1}^m \lambda_j =- \frac{a_{m-1}}{a_m}
\end{equation}
such that eq. (\ref{SumHj}) is proved. 
\vspace{0.5cm}

\noindent 9c) Gershgorin’s circle theorems by means of the matrix 
\begin{equation}
\mathbf{F}(\lambda )= \mathbf{E} - \lambda\, \mathbf{I}_m . 
\end{equation}
Let be a circle with the center point 
\begin{equation}
H_k=U_k + V_k\cdot i
\end{equation}
and the radius 
\begin{equation}
r_k= (n-1)\cdot \vert d_k\vert . 
\end{equation}
If the circle is separated from the remaining $n-1$ circles then we have to distinguish two cases
\vspace{0.5cm}

\noindent 9c1) The main value $H_k$ is real. Then also the included eigenvalue $\lambda_k$ is real and we have 
\begin{equation}
-r_k + H_k < \lambda_k < H_k + r_k .
\end{equation}

\noindent 9c2) For a complex eigenvalue
\begin{equation}
\lambda_k= u_k + v_k\cdot i
\end{equation}
we have the enclosures
\begin{equation}
-r_k + U_k < u_k < U_k + r_k
\end{equation}
and 
\begin{equation}
-r_k +V_k < v_k < V_k + r_k .
\end{equation} 
In the case of multiple eigenvalues or eigenvalue clusters we have simultaneous enclosures; cf. (\cite{Falk2004}, p.52).
\vspace{0.5cm}

\noindent 9d) Order reduction

\noindent 9d1) Scalar Polynomial. Separated a zero using Horner's scheme. 

\noindent 9d2) Matrix polynomial (\ref{PM}). Separated a cluster of $n$ eigenvalues en bloc [5].

\section{Numerical Examples}

\noindent
{\bf Example 1:} 

\noindent Following section \ref{exploration} an exploration is performed by means of the Pad\'e function (\ref{Pade}). The polynomial 
\begin{equation}
f(\lambda )= 4+12\lambda + 9 \lambda^2 - 4 \lambda^3 - 6 \lambda^4 + 0\cdot \lambda^5+\lambda^6
\end{equation}
is assumed with the zeros
\begin{equation}
\lambda_1=\lambda_2=2;\ \lambda_3=\lambda_4=\lambda_5=\lambda_6=-1.
\end{equation}
An exploration with $\delta = 0.3$ results in the pairs of values
\vspace{0.3cm}

\begin{equation}
\noindent\begin{tabular}{_l*{3}{^l}}
 &$\quad\quad\quad\quad\Lambda$               & $\quad\quad\quad\quad p(\Lambda )$ &      \\
\midrule
& $3.000000000000000e-01$ & $-5.261904761904761e-01$   & \rdelim\}{5}{0cm} \\
& $6.000000000000000e-01$ & $-9.333333333333332e-01$  & \\
& $9.000000000000000e-01$ & $-3.483333333333336e+00$  & \\
& $1.200000000000000e+00$ & $+1.466666666666667e+00$  & \\
& $1.500000000000000e+00$ & $+4.166666666666667e-01$  & \\
\midrule
$\lambda_1=$& $1.800000000000000e+00$ & $+1.166666666666663e-01$  & \rdelim\}{2}{0cm}[\normalfont change of sign] \\
$\lambda_2=$& $2.100000000000000e+00$ & $-4.696969696969665e-02$  &  \\
\end{tabular}
\end{equation}
\vspace{0.3cm}

\noindent It follows the regula falsi (\ref{delta2}) with a rounding to 
\begin{equation}
\lambda_3=2.01389
\end{equation}
and therefore with (\ref{Lambdapv}) 
\begin{equation}
\begin{tabular}{_l*{3}{^l}}
$j$ &$\quad\quad\quad\quad\Lambda_j$               & $\quad\quad\quad\quad p_1(\Lambda_j )$  &      \\
\midrule
$1$ & $2.013890000000000e+00$ & $-3.428770022238466e-03$   & \\
$2$ & $2.006984834339914e+00$ & $-1.735089829011837e-03$  & \\
$3$ & $2.003502535366870e+00$ & $-8.728302942998172e-04$  & \\
$4$ & $2.001753817659295e+00$ & $-4.377505086146021e-04$  & \\
$5$ & $2.000877548907494e+00$ & $-2.192108705613412e-04$  & \\
\end{tabular}
\end{equation}
There is no convergence for $\nu=1$.
\begin{equation}
\begin{tabular}{_l*{3}{^l}}
$j$ &$\quad\quad\quad\quad\Lambda_j$               & $\quad\quad\quad\quad p_2(\Lambda_j )$  &      \\
\midrule
$1$ & $2.013890000000000e+00$ & $-6.734450893113345e-03$   & \\
$2$ & $2.000327556690868e+00$ & $-1.636577281729843e-04$  & \\
$3$ & $2.000000187627338e+00$ & $-9.381362945114195e-08$  & \\
$4$ & $2.000000000000062e+00$ & $-3.090806074727201e-14$  & \\
$5$ & $2.000000000000000e+00$ & $0$  & \\
\end{tabular}
\end{equation}
It converges for $\nu=2$ and therefore we have
\begin{equation}
\lambda_1=\lambda_2=2.
\end{equation}
If $\lambda$ is replaced by $-\lambda$ we obtain the co-polynomial
\begin{equation}
f(-\lambda)= 4 - 12 \lambda + 9 \lambda^2 + 4 \lambda^3 - 6 \lambda^4 - 0\cdot \lambda^5 + \lambda^6 .
\end{equation}
An exploration results in
\begin{equation}
\noindent\begin{tabular}{_l*{3}{^l}}
 &$\quad\quad\quad\quad\Lambda$               & $\quad\quad\quad\quad p(\Lambda )$ &      \\
\midrule
& $3.000000000000000e-01$ & $-2.064102564102564e-01$   & \\
& $6.000000000000000e-01$ & $-1.083333333333336e-01$  & \\
& $9.000000000000000e-01$ & $-2.543859649125041e-02$  & \\
& $1.200000000000000e+00$ & $+4.848484848484354e-02$  & \\
& $1.500000000000000e+00$ & $+1.166666666666667e-01$  & \\
\end{tabular}
\end{equation}
It follows the regula falsi (\ref{delta2}) with a rounding to 
\begin{equation}
\lambda_3=1.00324 .
\end{equation}
With algorithm (\ref{Lambdapv}) we obtain no convergence for $\nu=1$, $\nu=2$, $\nu=3$ but
for $\nu=4$
\begin{equation}
\begin{tabular}{_l*{3}{^l}}
$j$ &$\quad\quad\quad\quad\Lambda_j$               & $\quad\quad\quad\quad p_4(\Lambda_j )$  &      \\
\midrule
$1$ & $1.003240000000000e+00$ & $-3.191729984318953e-03$   & \\
$2$ & $1.000037928810532e+00$ & $-3.792209825976858e-05$  & \\
$3$ & $1.000000005273932e+00$ & $-5.273931558410046e-09$  & \\
$4$ & $1.000000000000000e+00$ & $-2.343804163097548e-16$  & \\
\end{tabular}
\end{equation}
We have four times $+1$ of the co-polynomial and therefore
\begin{equation}
\lambda_1=\lambda_2=\lambda_3=\lambda_4=-1.
\end{equation}
With MATLAB the following zeros are calculated

\begin{eqnarray}\nonumber
\tilde{\lambda}_1 &=& + 2.000000000000001e+00+7.152216756864169e-09i  \nonumber\\
\tilde{\lambda}_2 &=& + 2.000000000000001e+00-7.152216756864169e-09i      \nonumber\\
\tilde{\lambda}_3 &=&  -1.000143391292847e+00      \\
\tilde{\lambda}_4 &=& -9.999999991419022e-01+1.433904397109860e-04i     \nonumber \\
\tilde{\lambda}_5 &=& -9.999999991419022e-01-1.433904397109860e-04i     \nonumber \\
\tilde{\lambda}_6 &=& -9.998566104233441e-01      \nonumber
\end{eqnarray}

\noindent
{\bf Example 2:} 

\noindent Let us consider a diagonal dominant matrix 

\begin{equation} %
\mathbf{F}(\lambda )=\pmatrix{   %
5 + 2 \lambda +3{\lambda}^2 &-1& 0& 0& 0 \cr
-1& 9+3\lambda +\lambda^2&-3&-2& 0 \cr
0& -3& 6+\lambda^2 & -2& 0 \cr
0&-2 & 2& 12+\lambda +\lambda^2& -5 -\lambda \cr
0& 0& 0 & -5-\lambda & 8+4\lambda+4\lambda^2
}
\end{equation}
with 
\begin{eqnarray}
\hbox{det} \mathbf{F}(\lambda ) = f(\lambda )&=& 1221 + 19366 \lambda + 33492 \lambda^2 + 28079 \lambda^3 + 23637 \lambda^4 \nonumber\\
&& + 11574 \lambda^5 + 5699 \lambda^6 + 1631 \lambda^7 + 489 \lambda^8 + 68 \lambda^9 + 12 \lambda^{10}
\end{eqnarray}
The quadratic equations (\ref{quadraticEq})
\begin{eqnarray}
5 + 2\lambda + 3 \lambda^2 &=& 0\nonumber \\
9+3\lambda +\lambda^2&=& 0\nonumber\\
6+\lambda^2&=& 0\\
12+\lambda +\lambda^2&=& 0\nonumber \\
8+4\lambda+4\lambda^2&=& 0\nonumber 
\end{eqnarray}
have the zeros
\begin{eqnarray}\label{zeros}
&& -3.333333333333334e-01+1.247219128924647e+00 i\nonumber\\
&& -3.333333333333334e-01-1.247219128924647e+00 i\nonumber\\
&& -1.500000000000000e+00+2.598076211353316e+00 i\nonumber\\ 
&& -1.500000000000000e+00-2.598076211353316e+00 i\nonumber\\
&& \quad 0\quad\quad\quad\quad\quad\quad\quad\quad\quad\quad\ \ +2.449489742783178e+00 i\\
&& \quad 0\quad\quad\quad\quad\quad\quad\quad\quad\quad\quad\ \ -2.449489742783178e+00 i\nonumber\\
&& -4.999999999999998e-01+3.427827300200522e+00 i\nonumber\\
&& -4.999999999999998e-01-3.427827300200522e+00 i\nonumber\\
&& -5.000000000000000e-01+1.322875655532295e+00 i\nonumber\\
&& -5.000000000000000e-01-1.322875655532295e+00 i\nonumber
\end{eqnarray}
We choose the third row as from (\ref{zeros}) as starting value and obtain the following results.

\noindent a) Algorithm (\ref{algor21})
\begin{equation}
\begin{tabular}{_l*{3}{^l}}
$j$ &$\quad\quad\quad\quad\Lambda_j$               & $\quad\quad\quad\quad p_4(\Lambda_j )$        \\
\midrule
$1$ & $-1.500000000000000e+00+2.598076211353316e+00 i$ & $+2.032105570683291e-01-7.395719545396148e-02 i$    \\
$2$ & $-1.296789442931671e+00+2.524119015899355e+00 i$ & $+1.773110756693907e-01-2.275434281511629e-02 i$   \\
$3$ & $-1.119478367262280e+00+2.501364673084238e+00 i$ & $+1.320427191450034e-01+5.627450679466858e-02 i$   \\
$4$ & $-9.874356481172767e-01+2.557639179878907e+00 i$ & $-3.007105982831948e-02+8.920521695652191e-02 i$ \\
$5$ & $-1.017506707945596e+00+2.646844396835429e+00 i$ & $+7.359853863202630e-04-2.137437648438716e-02 i$ \\
$6$ & $-1.016770722559276e+00+2.625470020351042e+00 i$ & $-9.748769601593532e-04-1.081681076639398e-03 i$\\
$7$ & $-1.017745599519435e+00+2.624388339274403e+00 i$ & $-5.137187053597404e-06+4.029601688399558e-06 i$\\
$8$ & $-1.017750736706489e+00+2.624392368876091e+00 i$ & $+1.136121179599617e-10-6.578228494562704e-11 i$\\
$9$ & $-1.017750736592877e+00+2.624392368810308e+00 i$ & $+2.009798662754188e-15-1.502768992132371e-15 i$
\end{tabular}
\end{equation}

\noindent b) Algorithm (\ref{Lambdapv})
\begin{equation}
\begin{tabular}{_l*{3}{^l}}
$j$ &$\quad\quad\quad\quad\Lambda_j$               & $\quad\quad\quad\quad p_4(\Lambda_j )$        \\
\midrule
$1$ & $-1.500000000000000e+00+2.598076211353316e+00 i$ & $-5.590550503850626e-02-5.178554874365029e-02 i$    \\
$2$ & $-1.281598940159485e+00+2.530507771744555e+00 i$ & $-3.206896765146202e-02-5.435962387848091e-02 i$   \\
$3$ & $-1.102941934511756e+00+2.519024236220835e+00 i$ & $+4.957773050493467e-03-4.558545263614854e-02 i$   \\
$4$ & $-9.935792103013821e-01+2.581791094008818e+00 i$ & $+2.032174795110959e-02+2.087166993263317e-03 i$ \\
$5$ & $-1.019159105737506e+00+2.632183836150752e+00 i$ & $-2.713265716958361e-03+4.879297147422435e-04 i$ \\
$6$ & $-1.017678176984104e+00+2.624544743975652e+00 i$ & $-4.116398540465768e-05+4.357842452009323e-05 i$\\
$7$ & $-1.017745599519435e+00+2.624388339274403e+00 i$ & $-5.137187053597404e-06+4.029601688399558e-06 i$\\
$8$ & $-1.017750658819505e+00+2.624392358442496e+00 i$ & $+1.342417115296137e-08+2.442886475494064e-08 i$\\
$9$ & $-1.017750736592890e+00+2.624392368810295e+00 i$ & $+2.991635783682046e-15-6.560209464970359e-15 i$
\end{tabular}
\end{equation}
Both algorithms converge quadratic and deliver almost identical results. 
In the same manner the remaining nine zeros will be calculated in parallel and independent from each other. 
\vspace{0.5cm}

\noindent
{\bf Example 3:} 

\noindent Evolution following section \ref{evolution}.

\noindent The polynomial denoted after Wilkinson
\begin{equation}\label{Wilkinson}
f(\lambda )= (1-\lambda)(2-\lambda)\cdots (9-\lambda)(10-\lambda)
\end{equation}
\noindent or in a decomposed form
\begin{eqnarray}
f(\lambda )&=& 3828800-10628640 \lambda +12753576 \lambda^2 - 8409500 \lambda^3 + 3416930 \lambda^4 - 902055 \lambda^5 + \\
& & + 157773 \lambda^6 - 18150 \lambda^7 + 1320 \lambda^8 - 55 \lambda^9 + \lambda^{10} 
\end{eqnarray}
has the zeros $1,2,\ldots ,10$.

\noindent Calculated zeros with MATLAB

\begin{eqnarray}
\tilde{\lambda}_1 &=&  1.000000000032865e+01  \nonumber\\
\tilde{\lambda}_2 &=&  8.999999998364443e+00 \nonumber\\
\tilde{\lambda}_3 &=&  8.000000003420013e+00 \nonumber\\
\tilde{\lambda}_4 &=&  6.999999996085851e+00  \nonumber\\
\tilde{\lambda}_5 &=&  6.000000002669752e+00  \\
\tilde{\lambda}_6 &=&  4.999999998898655e+00 \nonumber \\
\tilde{\lambda}_7 &=&  4.000000000263102e+00 \nonumber \\
\tilde{\lambda}_8 &=&  2.999999999968169e+00 \nonumber \\
\tilde{\lambda}_9 &=&  2.000000000001345e+00  \nonumber\\
\tilde{\lambda}_{10} &=& 1.000000000000000e+00  \nonumber
\end{eqnarray}
We use these values as interpolation values and obtain the following list
\begin{equation}\label{list_Ex3_1}
L_{10}=\pmatrix{\hbox{Interpolation Values}&\hbox{Defects}&\hbox{Main Values}\cr 
\midrule\cr
1.000000000032865e+01&+3.727125322099494e-10& 9.999999999955941e+00 \cr
8.999999998364443e+00&-1.720094133611112e-09& 9.000000000084537e+00 \cr
8.000000003420013e+00&+3.167697847832428e-09& 8.000000000252316e+00 \cr
6.999999996085851e+00&-4.044785689387324e-09& 7.000000000130637e+00 \cr
6.000000002669752e+00&+2.348194056725277e-09& 6.000000000321559e+00 \cr
4.999999998898655e+00&-1.197945997413012e-09& 5.000000000096601e+00 \cr
4.000000000263102e+00&+2.777798930869391e-10& 3.999999999985322e+00 \cr
2.999999999968169e+00&-3.340011018696152e-11& 3.000000000001569e+00 \cr
2.000000000001345e+00&+1.212659602373197e-12& 2.000000000000132e+00 \cr
1.000000000000000e+00&+0.000000000000000e+00& 1.000000000000000e+00
}
\end{equation}
It follows two evolutions
\begin{equation}\label{list_Ex3_2}
L_{10}=\pmatrix{\hbox{Interpolation Values}&\hbox{Defects}&\hbox{Main Values}\cr 
\midrule\cr
9.999999999955941e+00&-5.567454977060759e-11& 1.000000000001162e+01 \cr
9.000000000084537e+00&-3.564064318273009e-11& 9.000000000120178e+00 \cr
8.000000000252316e+00&-8.075158037691484e-11& 8.000000000333067e+00 \cr
7.000000000130637e+00&+2.082540757136542e-10& 6.999999999922383e+00 \cr
6.000000000321559e+00&+2.186021042651543e-10& 6.000000000102957e+00 \cr
5.000000000096601e+00&+6.758556180612028e-11& 5.000000000020015e+00 \cr
3.999999999985322e+00&-1.347399557668777e-11& 3.999999999998797e+00 \cr
3.000000000001569e+00&+4.296279733101633e-12& 2.999999999997273e+00 \cr
2.000000000000132e+00&+6.929483440822838e-14& 2.000000000000063e+00 \cr
1.000000000000000e+00&+0.000000000000000e+00& 1.000000000000000e+00
}
\end{equation}
\begin{equation}\label{list_Ex3_3}
L_{10}=\pmatrix{\hbox{Interpolation Values}&\hbox{Defects}&\hbox{Main Values}\cr 
\midrule\cr
1.000000000001162e+01&+1.774461056701600e-11& 9.999999999993872e+00 \cr
9.000000000020178e+00&-7.047284661419908e-11& 9.000000000190651e+00 \cr
8.000000000333067e+00&-7.797978696197630e-11& 8.000000000411047e+00 \cr
6.999999999922383e+00&-3.380356010780448e-10& 7.000000000260418e+00 \cr
6.000000000102957e+00&+1.409918897505596e-10& 5.999999999961965e+00 \cr
5.000000000029015e+00&+7.793359042933322e-11& 4.999999999951082e+00 \cr
3.999999999998797e+00&-5.928558055009481e-12& 4.000000000004725e+00 \cr
2.999999999997273e+00&-2.910383045325759e-12& 3.000000000000184e+00 \cr
2.000000000000063e+00&+4.619655627634508e-14& 2.000000000000016e+00 \cr
1.000000000000000e+00&+0.000000000000000e+00& 1.000000000000000e+00
}
\end{equation}
Using the sum control (\ref{SumHj}) we have ($a_{10}=1$)
\begin{equation}
\hbox{Desired value:} \sum_{j=1}^{10} H_j = -\frac{-55}{a_{10}}=55.\quad \hbox{Actual value:}\ 49.466
\end{equation}
\vspace{0.5cm}

\noindent
{\bf Example 4:} 

\noindent Correction of multiple eigenvalues after (\ref{Lambdapv}).

\noindent Let be a matrix pencil
\begin{equation}
{\mathbf{F}(\lambda)=\mathbf{A}-\lambda \mathbf{B}}
\end{equation}
with
\begin{equation}
\mathbf{A}= \pmatrix{-1 & 0& 1 & 0 &0 \cr
0& 0& 0& 1 & 0 \cr
1& 0& 0& 0 & 1 \cr
0& 1& 0& 0 & 0 \cr
0& 0& 0& 1 & -1 \cr}; \mathbf{B}=\mathbf{I}_5 .
\end{equation}
Its eigenvalues are 
\begin{equation}
\lambda_1=-2;\ \lambda_2 = \lambda_3 = -1;\ \lambda_4 = \lambda_5 = 1.
\end{equation}
\begin{figure}[h]
\centering
\includegraphics[width=0.5\linewidth]{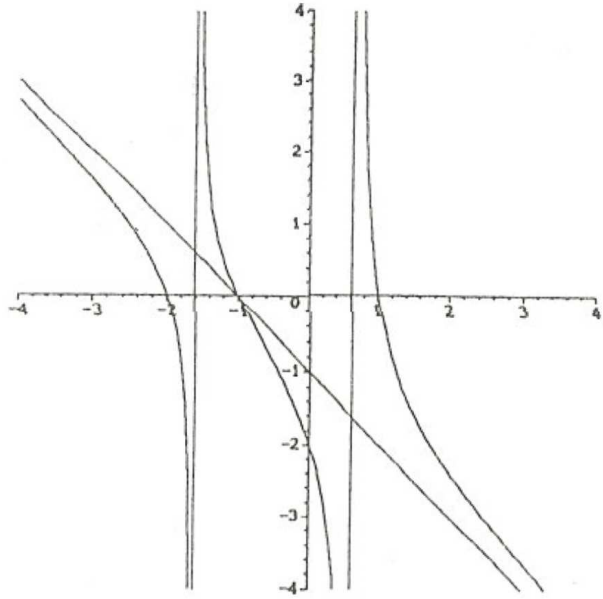}
\caption[]{The Pad\'e function of the matrix pencil of example 4}
\label{fig:Falk}
\end{figure}

MATLAB calculates the following approximated zeros 
\begin{eqnarray}
\tilde{\lambda}_1 &=&  -1.999999999999996e+00  \nonumber\\
\tilde{\lambda}_2 &=&  +1.000000000000000e+00+7.768125062636118e-09i \nonumber\\
\tilde{\lambda}_3 &=&  +1.000000000000000e+00-7.768125062636118e-09i \\
\tilde{\lambda}_4 &=&  -1.000000009896685e+00  \nonumber\\
\tilde{\lambda}_5 &=&  -9.999999901033162e-01  \nonumber\\
\end{eqnarray}
It follows the corrections after (\ref{Lambdapv}). We start with $\tilde{\lambda}_5$.

\begin{equation}
\begin{tabular}{_l*{3}{^l}}
$j$ &$\quad\quad\quad\quad\Lambda_j$               & $\quad\quad\quad\quad p_1(\Lambda_j )$  &      \\
\midrule
$1$ & $-9.999999901033162e-01$ & $-1.121813169708485e-08$   & \\
$2$ & $-1.000000001321448e+00$ & $0$  & \\
\end{tabular}
\end{equation}
No convergence for $\nu =1$.
\begin{equation}
\begin{tabular}{_l*{3}{^l}}
$j$ &$\quad\quad\quad\quad\Lambda_j$               & $\quad\quad\quad\quad p_2(\Lambda_j )$  &      \\
\midrule
$1$ & $-9.999999901033162e-01$ & $-9.896683722532278e-09$   & \\
$2$ & $-9.999999999999999e-01$ & $-1.665334536937735e-16$  & \\
\end{tabular}
\end{equation}
Convergence for $\nu =2$. Therefore, we have $\lambda_4=\lambda_5=1$. 

\noindent We start with $\tilde{\lambda}_1$:
 
 \begin{equation}
 \begin{tabular}{_l*{3}{^l}}
 $j$ &$\quad\quad\quad\quad\Lambda_j$               & $\quad\quad\quad\quad p_1(\Lambda_j )$  &      \\
 \midrule
 $1$ & $-1.999999999999996e+00$ & $-3.552713678800562e-15$   & \\
 $2$ & $-2.000000000000000e+00$ & $0$  & \\
 \end{tabular}
 \end{equation}
Therefore we have $\lambda_1=-2$.

\noindent Now, we start with $\tilde{\lambda}_2$:

 \begin{equation}
 \small
 \begin{tabular}{_l*{3}{^l}}
 $j$ &$\quad\quad\quad\quad\Lambda_j$               & $\quad\quad\quad\quad p_1(\Lambda_j )$       \\
 \midrule
 $1$ & $1.000000000000000e+00+7.768125062636118e-09 i$ & $-1.491714852512857e-16 -4.764011116681661e-09 i$  \\
 $2$ & $1.000000000000000e+00+3.004113945954457e-09 i$ & $-2.065119601239563e-16 -1.334040524532894e-23 i$   \\
 \end{tabular}
 \end{equation}
 No convergence for $\nu =1$.
 
  \begin{equation}
  \small
  \begin{tabular}{_l*{3}{^l}}
  $j$ &$\quad\quad\quad\quad\Lambda_j$               & $\quad\quad\quad\quad p_2(\Lambda_j )$       \\
  \midrule
  $1$ & $1.000000000000000e+00+7.768125062636118e-09 i$ & $-5.744419753425677e-16 -7.768125062636107e-09 i$  \\
  $2$ & $1.000000000000000e+00+3.004113945954457e-09 i$ & $+1.295260195396017e-16 -1.158052857574239e-23 i$   \\
  \end{tabular}
  \end{equation}
Convergence for $\nu =2$.

\noindent Therefore, we have $\lambda_2=\lambda_3=-1$.
\vspace{0.5cm}

%

\noindent
{\bf Example 5:}

\noindent Multiple complex zeros: The polynomial 
\begin{equation}
f(\lambda )= \left( 1 + \lambda + \lambda^2   \right)^3\cdot \left( 1  + \lambda^2   \right)^2\cdot 6
\end{equation}
can be decomposed into
\begin{equation}\label{PolynEx5}
f(\lambda )= 6+ 18 \lambda + 48 \lambda^2 + 78 \lambda^3 + 114 \lambda^6 + 78 \lambda^7 + 48 \lambda^8 + 18 \lambda^9 + 6 \lambda^{10}.
\end{equation}
The zeros are 
\begin{eqnarray}
\lambda_1 &=&    0+i\nonumber\\
\lambda_2 &=&    0-i\nonumber\\
\lambda_3 &=&   0+i\nonumber\\
\lambda_4 &=&   0-i \nonumber\\
\lambda_5 &=&  -0,5 + \sqrt{0,75} i  \label{PolynEx5_zero}\\
\lambda_6 &=&  -0,5 - \sqrt{0,75} i \nonumber \\
\lambda_7 &=&   -0,5 + \sqrt{0,75} i\nonumber \\
\lambda_8 &=&   -0,5 - \sqrt{0,75} i\nonumber \\
\lambda_9 &=&    -0,5 + \sqrt{0,75} i\nonumber\\
\lambda_{10} &=&    -0,5 - \sqrt{0,75} i\nonumber
\end{eqnarray}
with 
\begin{equation}
\sqrt{0.75}=8.660254037844386 e-01 .
\end{equation}
MATLAB calculates the following approximated zeros
\begin{eqnarray}
\tilde{\lambda}_1 &=& +2.103940549558203e-08+1.000000028920264e+00 i    \nonumber\\
\tilde{\lambda}_2 &=& +2.103940549558203e-08-1.000000028920264e+00 i   \nonumber\\
\tilde{\lambda}_3 &=& -2.103939766850971e-08+9.999999710797240e-01 i  \nonumber\\
\tilde{\lambda}_4 &=& -2.103939766850971e-08-9.999999710797240e-01 i   \nonumber\\
\tilde{\lambda}_5 &=& -5.000094136551562e-01+8.660276783463672e-01 i   \\
\tilde{\lambda}_6 &=& -5.000094136551562e-01-8.660276783463672e-01 i  \nonumber \\
\tilde{\lambda}_7 &=& -4.999933232335635e-01+8.660324192348879e-01 i  \nonumber \\
\tilde{\lambda}_8 &=& -4.999933232335635e-01-8.660324192348879e-01 i  \nonumber \\
\tilde{\lambda}_9 &=&  -4.999972631112927e-01+8.660161137720683e-01 i   \nonumber\\
\tilde{\lambda}_{10}&=&-4.999972631112927e-01-8.660161137720683e-01 i    \nonumber
\end{eqnarray}

\noindent Correction of $\tilde{\lambda}_5$ after (\ref{Lambdapv}).

\begin{equation}
 \small
 \begin{tabular}{_l*{3}{^l}}
 $j$ &$\quad\quad\quad\quad\Lambda_j$               & $\quad\quad\quad\quad p_1(\Lambda_j )$       \\
 \midrule
 $1$ & $-5.000094136551562e-01+8.660276783463672e-01 i$ & $+2.590606103860521e-06 -4.667337884496285e-07 i$  \\
 $2$ & $-5.000068230490523e-01+8.660272116125788e-01 i$ & $+5.061285460948595e-06 -9.728311147477021e-07 i$   \\
 $3$ & $-5.000017617635913e-01+8.660262387814640e-01 i$ & $-3.494617894702583e-05 +1.077378144656732e-05 i$  \\
 $4$ & $-5.000367079425384e-01+8.660370125629105e-01 i$ & $+1.215888288037037e-05 -3.935532380990325e-06 i$   \\
 $5$ & $-5.000245490596580e-01+8.660330770305296e-01 i$ & $+2.065119601239563e-06 -2.552780638469371e-06 i$   \\
 \end{tabular}
 \end{equation}
 No convergence for $\nu =1$.
 
  \begin{equation}
  \small
  \begin{tabular}{_l*{3}{^l}}
  $j$ &$\quad\quad\quad\quad\Lambda_j$               & $\quad\quad\quad\quad p_2(\Lambda_j )$       \\
  \midrule
 $1$ & $-5.000094136551562e-01+8.660276783463672e-01 i$ & $4.706686932467359e-06 -1.137278373809355e-06 i$  \\
 $2$ & $-5.000047069682237e-01+8.660265410679934e-01 i$ & $2.353445844274929e-06 -5.686590462229913e-07 i$   \\
 $3$ & $-5.000023535223794e-01+8.660259724089472e-01 i$ & $1.176760802740233e-06 -2.842232846214180e-07 i$  \\
 $4$ & $-5.000011767615767e-01+8.660256881856626e-01 i$ & $5.883348476664755e-07 -1.420315348699738e-07 i$   \\
 $5$ & $-5.000005884267291e-01+8.660255461541277e-01 i$ & $2.940921445040524e-07 -7.102102632529488e-08 i$   \\
  \end{tabular}
  \end{equation}
No convergence for $\nu =2$.

 \begin{equation}
  \small
  \begin{tabular}{_l*{3}{^l}}
  $j$ &$\quad\quad\quad\quad\Lambda_j$               & $\quad\quad\quad\quad p_3(\Lambda_j )$       \\
  \midrule
  $1$ & $-5.000000008552082e-01+8.660254038125426e-01 i$ & $+8.552084998311037e-10-2.810326544530710e-11 i$  \\
  $2$ & $-4.999999999999997e-01+8.660254037844393e-09 i$ & $-1.553244317250170e-15-2.433903049193858e-15 i$   \\
  \end{tabular}
  \end{equation}
Convergence for $\nu =3$ such that we have a zero $\lambda_5$ of (\ref{PolynEx5_zero}) with the multiplicity 3. The polynomial (\ref{PolynEx5}) is (accidental) hermitian but of even order $m=10$ and therefore $-1$ is no zero.
\vspace{0.5cm}

\noindent
{\bf Example 6:}

\noindent The reduced eigenvalues equation (\ref{redEV}) with 
\begin{equation}
p_E(\lambda)= \frac{S_1 (\lambda)-1}{-S_2 (\lambda)}
\end{equation}

\noindent Wilkinson polynomial (\ref{Wilkinson}):
\begin{equation}\label{list_Ex3_3}
L_{10}=\pmatrix{\hbox{Interpolation Values}&\hbox{Defects}&\hbox{Main Values}\cr 
\midrule\cr
1.000100000001162e+01&9.988178397837826e-05& 1.000000118216022e+00 \cr
2.000200000020178e+00&1.997514392797159e-04& 2.000000248560720e+00 \cr
3.000300000333067e+00&2.996318462367703e-04& 3.000000368153763e+00 \cr
4.000400000000000e+00&3.995415371807089e-04& 4.000000458462819e+00 \cr
5.000500000000000e+00&4.995001173401799e-04& 5.000000499882660e+00 \cr
6.000599999999999e+00&5.995318981743213e-04& 6.000000468101825e+00 \cr
7.000699999999999e+00&6.996716090711247e-04& 7.000000328390928e+00 \cr
8.000800000000000e+00&7.999787558876377e-04& 8.000000021244112e+00 \cr
9.000900000000063e+00&9.005807841077982e-04& 8.999999419215891e+00 \cr
1.000100000000000e+01&1.001930111657591e-03& 9.999998069888342e+00
}
\end{equation}

\noindent We start with the main value $H_3$ an obtain
\begin{equation}
  \small
  \begin{tabular}{_l*{3}{^l}}
  $j$ &$\quad\quad\quad\quad\Lambda_j$               & $\quad\quad\quad\quad p_2(\Lambda_j )$       \\
  \midrule
  $1$ & $3.000000368155010e+00$ & $-3.677030854107595e-07$  \\
  $2$ & $3.000000000451924e+00$ & $-4.517926727218366e-10$   \\
  $3$ & $3.000000000000131e+00$ & $-5.277591961897212e-16$   \\
  \end{tabular}
  \end{equation}
In comparison with the algorithm (\ref{algor21})
\begin{equation}
  \small
  \begin{tabular}{_l*{3}{^l}}
  $j$ &$\quad\quad\quad\quad\Lambda_j$               & $\quad\quad\quad\quad p_2(\Lambda_j )$       \\
  \midrule
  $1$ & $3.000000368155010e+00$ & $-3.681554823412160e-07$  \\
  $2$ & $2.999999999999527e+00$ & $+2.448417482860382e-12$   \\
  $3$ & $3.000000000001975e+00$ & $-8.777345693336323e-13$   \\
  \end{tabular}
  \end{equation}
\vspace{0.5cm}
 
 \noindent
 {\bf Example 7:}
 
 \noindent Singular leading matrix.
 \vspace{0.2cm}
 
 \noindent We assume a polynomial matrix
 \begin{equation}
 \mathbf{F}(\lambda) = \mathbf{A}_0 + \mathbf{A}_1 \lambda + \mathbf{A}_2 \lambda^2 + \mathbf{A}_3 \lambda^3+\mathbf{A}_4 \lambda^4
 \end{equation}
 with the coefficient matrices
 \begin{eqnarray}
 \mathbf{A}_0 &=&  \pmatrix{1&0\cr 0&1}, \nonumber\\
 \mathbf{A}_1 &=&  \pmatrix{1&1\cr 1&1}, \nonumber\\
 \mathbf{A}_2 &=&  \pmatrix{2&1\cr 0&1}, \\
 \mathbf{A}_3 &=&  \pmatrix{0&0\cr 0&0}, \nonumber\\
 \mathbf{A}_4 &=&  \pmatrix{0&1\cr 0&0},\nonumber
 \end{eqnarray}
 where the leading matrix $\mathbf{A}_4$ is singular such that we have fewer than $m=\rho\cdot n$ eigenvalues. 
 \vspace{0.2cm}
 
 \noindent We have to distinguish two approaches: 
 \vspace{0.2cm}
 
 \noindent a) using the matrix 
 \begin{equation}\label{Ex7matrix}
 \mathbf{F}(\lambda)=
\left( \begin{array}{cc}
1+\lambda + 2 \lambda^2 & \lambda + 2 \lambda^2+\lambda^4\\
 \lambda & 1+\lambda + 2 \lambda^2
\end{array} \right)
 \end{equation}
 \vspace{0.2cm}
 
 \noindent using the characteristic polynomial 
 \begin{equation}\label{Ex7polyn}
 \hbox{det}\,\mathbf{F}(\lambda)=f(\lambda )= 1 + 2 \lambda + 3 \lambda^2 + 2\lambda^3 + 2 \lambda^4 - \lambda^5
 \end{equation}
 with the degree 5; therefore, we have only 5 zeros and accordingly 5 eigenvalues.
 \vspace{0.2cm}
 
\noindent We start with the exploration using the Pad\'e function (\ref{Pade}) and choose $\delta =0.1$:
\begin{equation}
\noindent\begin{tabular}{_l*{3}{^l}}
 &$\quad\quad\quad\quad\Lambda_j$               & $\quad\quad\quad\quad p(\Lambda_j )$ &      \\
\midrule
& $0.0$ & $-2.064102564102564e-01$   & \rdelim\}{30}{0cm}\\
& $0.1$ & $-2.064102564102564e-01$  & \\
& $0.2$ & $-2.064102564102564e-01$  & \\
& $0.3$ & $-2.064102564102564e-01$  & \\
& $0.4$ & $-2.064102564102564e-01$  & \\
& $0.5$ & $-2.064102564102564e-01$   & \\
& $0.6$ & $-2.064102564102564e-01$  & \\
& $0.7$ & $-2.064102564102564e-01$  & \\
& $0.8$ & $-2.064102564102564e-01$  & \\
& $0.9$ & $-2.064102564102564e-01$  & \\
& $1.0$ & $-2.064102564102564e-01$   & \\
& $1.1$ & $-2.064102564102564e-01$  & \\
& $1.2$ & $-2.064102564102564e-01$  & \\
& $1.3$ & $-2.064102564102564e-01$  & \\
& $1.4$ & $-2.064102564102564e-01$  & \\
& $1.5$ & $-2.064102564102564e-01$   & \\
& $1.6$ & $-2.064102564102564e-01$  & \\
& $1.7$ & $-2.064102564102564e-01$  & \\
& $1.8$ & $-2.064102564102564e-01$  & \\
& $1.9$ & $-2.064102564102564e-01$  & \\
& $2.0$ & $-2.064102564102564e-01$   & \\
& $2.1$ & $-2.064102564102564e-01$  & \\
& $2.2$ & $-2.064102564102564e-01$  & \\
& $2.3$ & $-2.064102564102564e-01$  & \\
& $2.4$ & $-2.064102564102564e-01$  & \\
& $2.5$ & $-2.064102564102564e-01$   & \\
& $2.6$ & $-2.064102564102564e-01$  & \\
& $2.7$ & $-2.064102564102564e-01$  & \\
& $2.8$ & $-2.064102564102564e-01$  & \\
& $2.9$ & $-2.064102564102564e-01$  & \\
 $\lambda_1=$& $3.0$ & $-2.064102564102564e-01$  & \rdelim\}{2}{0cm}[\normalfont change of sign] \\
 $\lambda_2=$& $3.1$ & $-2.064102564102564e-01$  &\\
\end{tabular}
\end{equation}
and hence with the regula falsi 
\begin{equation}
\lambda_3=3.05965871206409e+00 
\end{equation}
and furthermore after (\ref{psigma}) with $\sigma =5$
 \begin{equation}
 \lambda_4=3.056811621817845e+00 .
 \end{equation}
 \vspace{0.2cm}
 It follows the Pad\'e algorithm (\ref{algor21})
 \begin{equation}\label{Ex7lambda3}
   \small
   \begin{tabular}{_l*{3}{^l}}
   $j$ &$\quad\quad\quad\quad\Lambda_j$               & $\quad\quad\quad\quad p(\Lambda_j )$       \\
   \midrule
   $1$ & $3.056811621817845e+00$ & $-2.231403025508382e-06$  \\
   $2$ & $3.056809390414819e+00$ & $-5.754993999073329e-12$   \\
   $3$ & $3.056809390409065e+00$ & $-7.589857143243228e-17$   \\
   \end{tabular}
   \end{equation}
 
 \noindent a) MATLAB calculates the eigenvalues for the matrix (\ref{Ex7matrix})
 \begin{eqnarray}
 \tilde{\lambda}_1 &=& +3.056809390409061e+00    \nonumber\\
 \tilde{\lambda}_2 &=& -2.103940549558203e-08-1.000000028920264e+00 i   \nonumber\\
 \tilde{\lambda}_3 &=& -2.103939766850971e-08+9.999999710797240e-01 i  \\
 \tilde{\lambda}_4 &=& -2.103939766850971e-08-9.999999710797240e-01 i   \nonumber\\
 \tilde{\lambda}_5 &=& -5.000094136551562e-01+8.660276783463672e-01 i  \nonumber 
 \end{eqnarray}
 
 \noindent b) MATLAB calculates the zeros for the polynomial (\ref{Ex7polyn})
  \begin{eqnarray}
 \tilde{\lambda}_1 &=& +3.056809390409070e+00   \nonumber\\
 \tilde{\lambda}_2 &=& -2.103940549558203e-08-1.000000028920264e+00 i   \nonumber\\
 \tilde{\lambda}_3 &=& -2.103939766850971e-08+9.999999710797240e-01 i  \\
 \tilde{\lambda}_4 &=& -2.103939766850971e-08-9.999999710797240e-01 i   \nonumber\\
 \tilde{\lambda}_5 &=& -5.000094136551562e-01+8.660276783463672e-01 i  \nonumber 
  \end{eqnarray}
  Both MATLAB results as well as $\Lambda_3$ in (\ref{Ex7lambda3}) are comparable with respect to the accuracy
  \vspace{0.5cm}
  
 \noindent
 {\bf Example 8:}
 
 \begin{equation}
 f(\lambda )= (\lambda - 1)(\lambda -2)(\lambda -3)(\lambda -4)(\lambda -5)\cdot 3
 \end{equation}
 or in a decomposed form
 \begin{equation}
 f(\lambda )= -360 + 822 \lambda - 675 \lambda ^2 + 255 \lambda^3 - 45 \lambda^4 + 3 \lambda^5 .
 \end{equation}
 The exploration with $\delta =0.3$ delivers the pairs of values
 \begin{equation}\label{Ex7lambda3}
    \small
    \begin{tabular}{_l*{3}{^l}}
    $$ &$\ \Lambda_j$               & $\quad\quad\quad\quad p(\Lambda_j )$ &      \\
    \midrule
                & $0.0$ & $+4.379562043795621e-01$ &\rdelim\}{3}{0cm}  \\
                & $0.3$ & $+3.484061594869381e-01$  &  \\
                & $0.6$ & $+2.408279034112688e-01$ &   \\
  $\lambda_1=$ &  $0.9$ & $+8.366965417990657e-02$ & \rdelim\}{2}{0cm}[\normalfont change of sign]  \\
 $\lambda_2=$ &  $1.2$ & $-3.884787018255549e-01$  &  \\
    \end{tabular}
    \end{equation}
With the regula falsi algorithm the following value can be calculated
\begin{equation}\label{ExRegulaFalsi}
\lambda_3= 9.531631550437540e-01 .
\end{equation}
Now, we use Halley's algorithm after (\ref{Halley}).
\begin{equation}\label{Ex7lambda3}
   \small
   \begin{tabular}{_l*{3}{^l}}
   $j$ &$\quad\quad\quad\quad\Lambda_j$               & $\quad\quad\quad\quad h(\Lambda_j )$       \\
   \midrule
   $1$ & $9.000000000000000e-01$ & $+1.180808950230746e-01$  \\
   $2$ & $1.018080895023075e+00$ & $-1.738293306349123e-02$   \\
   $3$ & $1.000697961959583e+00$ & $-6.969460936557861e-04$   \\
   $4$ & $1.000001015865928e+00$ & $-1.015863777488736e-06$  \\
   $5$ & $1.000000000002150e+00$ & $-2.150576013567233e-12$   \\
   $6$ & $9.999999999999994e-01$ & $-7.894919286223337e-16$   \\
   $7$ & $1.000000000000000e+00$ & $+0.000000000000000e+00$   \\
   \end{tabular}
   \end{equation}
It is known from the theory \cite{MathisTS}: the convergence of Halley's algorithm is cubic for simple zeros and quadratic for multiple zeros. However, at least in this example a cubic convergence cannot be observed. 
\vspace{0.2cm}

\noindent In comparison: the accelerated regula falsi following (\ref{ExRegulaFalsi}) leads in five steps to the nearly exact solution
\begin{equation}\label{Ex7lambda3}
   \small
   \begin{tabular}{_l*{3}{^l}}
   $j$ &$\quad\quad\quad\quad\lambda_j$               & $\quad\quad\quad\quad p(\lambda_j )$       \\
   \midrule
   $4$ & $9.985736255069474e-01$ & $+1.422152519655597e-03$  \\
   $5$ & $9.999780098758768e-01$ & $+2.198911675590514e-05$   \\
   $6$ & $1.000000002011310e+00$ & $-2.011309654213950e-09$   \\
   $7$ & $1.000000000000000e+00$ & $+0.000000000000000e+00$  \\
  
   \end{tabular}
   \end{equation}
\vspace{1cm}

\noindent {\bf Acknowlegdement:} I would like to thank M. Sc. Jonas Deni\ss en, Berlin, and Prof. Dr.-Ing. Wolfgang Mathis, Leibniz University of Hannover, for the successful cooperation.


\begin{thebibliography}{99}
 
\bibitem{CarstensenStein1987} %
Carstensen, C.; E. Stein (1987): Analysis and Calculation of Falk's ECP-transformation and familiar problems (in German). Intern. Series Numerical Mathem., vol. 83, pp. 47-61
\bibitem{CarstensenStein1989} %
Carstensen, C.; E. Stein (1989): About Falk's ECP-transformation and generalizations (in German). ZAMM, vol. 69, no. 11, pp. 375-391
\bibitem{Falk2004} %
Falk, S. (2004): The eigenvalue algorithm ECP for polynomial matrices (in German).
Abhandl. Braunschw. Wiss. Gesellsch., vol. 54
\bibitem{Falk2011} %
Falk, S. (2011): The accelerated Ritz iteration for polynomials and polynomial matrices (in German). Abhandl. Braunschw. Wiss. Gesellsch., vol. 58
\bibitem{Falk1994} %
Falk, S. (1994): The reduction method for polynomial matrices (in German). ZAMM, vol. 74, no. 1, pp. 3-15
\bibitem{MathisTS} %
Mathis, W.; T. Thiessen; S. Falk (2010): Iteration methods for zero problems (in German). Institute of Theoretical Electrical Engineering, Leibniz University of Hannover
\bibitem{ShawTraub} %
Shaw, M.; J.F. Traub (1974): In the number of multiplications for the evaluation of a polynomial and some ot its derivations. J. Assoc. Comp. mech., vol. 21, pp. 61-166
\bibitem{Stewart1968} %
Stewart, G. W. (1968): Some topics in numerical analysis. I Lehmer's method for finding the zeros of a polynomial. Techn. report, Oak Ridge
\bibitem{ZurmuehlFalk1} %
Zurm\"uhl, R.; S. Falk (1992): Matrices and its applications 1 (in German). 6. edition, Springer-Verlag, Berlin
\bibitem{ZurmuehlFalk2} %
Zurm\"uhl, R.; S. Falk (1986): Matrices and its applications 2 (in German). 5. edition, Springer-Verlag, Berlin





\end{thebibliography}
\end{document}